\documentclass{amsart}

\newtheorem{thm}{Theorem}

\theoremstyle{remark}
\newtheorem{rem}[thm]{Remark}

\title{Finite sums of commutators}
\author{Ciprian Pop}
\address{I.M.A.R., CP 1--764, Bucharest, Romania}
\curraddr{Department of Mathematics, Texas A\&M University, College
  Station, Texas 77843--3368} 
\email{cpop@math.tamu.edu}

\begin{document}
\begin{abstract}
  We show that elements of unital $C^*$-algebras without tracial states are
  finite sums of commutators. Moreover, the number of commutators involved
  is bounded, depending only on the given $C^*$-algebra.
\end{abstract}
\maketitle

\begin{section}{Introduction}
  It was shown in \cite{fack_harpe} that in finite von Neumann algebras,
  elements with central trace zero are sums of at most 10 commutators.
  The $C^*$-algebra case was considered in \cite{fack}. The main result
  there states that if the unit of a $C^*$-algebra $A$ is properly infinite
  (i.e.  there exists two orthogonal projections $p,q\in A$ equivalent to
  $1$), then any hermitian element is a sum of at most five self-adjoint
  commutators. In this paper we consider the more general case of unital
  $C^*$-algebras $A$ without tracial states and improve the previous result
  of T. Fack. Note that if the unit of $A$ is properly infinite, then $A$
  has no tracial states. The converse is known to be false, at least when
  $A$ is non-simple (see \cite{rordam} for further details). $C^*$-algebras
  without tracial states have several nice characterizations, such as
  \cite{haagerup}. This paper
  contains also another simple proof of the latter result of \cite[Lemma
  1]{haagerup}.
\end{section}

\begin{section}{The result}
  Given $a,b\in A$, their commutator is $[a,b]=ab-ba$. A self-adjoint
  commutator is just a commutator of the form $[a^*,a]=a^*a-aa^*$.

  \begin{thm}
    Let $A$ be a unital $C^*$-algebra. Then the following are
    equivalent:
    \begin{enumerate}
    \item $A$ has no tracial states.
    \item There exist an integer $n\ge2$, and elements
      $b_1,b_2,\ldots,b_n\in A$ such that
      \begin{eqnarray*}
        \sum_ib_i^*b_i&=&1\\
        \|\sum_ib_ib_i^*\|&<&1
      \end{eqnarray*}
    \item There exist an integer $n\ge2$ such that any element of $A$ can
      be expressed as a sum of $n$ commutators and any positive element can
      be expressed as a sum of at most $n$ self-adjoint commutators.
    \end{enumerate}
  \end{thm}
  
  \begin{rem}
    The equivalence of $1$ and $2$ is just \cite[Lemma 1]{haagerup}.
    As mentioned, we give here a new simple proof.
  \end{rem}
  \begin{rem}
    The integer $n$ appearing in both $2$ and $3$ above represents the
    same number. Thus, if the unit of $A$ is properly infinite, there exists
    two isometries $v_1,v_2\in A$ with orthogonal ranges. Let
    $b_i=(1/\sqrt{2})v_i$ for $i=1,2$. Then $b_1^*b_1+b_2^*b_2=1$ and
    $b_1b_1^*+b_2b_2^*=1/2$, thus the property from $2$ is achieved
    with $n=2$. Therefore positive elements are sums of $2$
    commutators and self-adjoint elements are sums of 4 commutators.
  \end{rem}

  \begin{proof}
    The implication $(3)\Rightarrow(1)$ is trivial.

    $(1)\Rightarrow(2)$. Consider 
    \begin{equation*}
      R=\{\sum_{i=1}^s(a^*_ia_i-a_ia_i^*)\ ;\ s\ge1,a_i\in A\}
    \end{equation*}
    the set of \textbf{finite} sums of self-adjoint commutators of
    $A$. Note that $R\subset A_{sa}$ is a \textbf{real vector subspace} of
    $A_{sa}$. Put $\delta=\mbox{dist}(1,R)$. 

    We show that $\delta<1$. Suppose the contrary, i.e. $\delta=1$. This is
    equivalent to 
    \begin{equation*}
      \|t+x\|\ge|t|,\quad\forall x\in R.
    \end{equation*}
    It follows that $\varphi(t+x)=t$ is a real bounded functional on
    $\mathbb{R}1+R$ of norm $1$. By the Hahn--Banach theorem it can be
    extended to a norm-$1$ functional on $A_{sa}$ and furthermore to a
    bounded \textbf{complex} functional on $A$, denoted also by $\varphi$.
    Observe that $\varphi$ is necessarily a tracial state on $A$, which
    contradicts our hypothesis.
    
    Because $\delta<1$, there exist some elements $a_1,a_2,\ldots,a_m\in A$
    such that $t_0=\|1-\sum_{i=1}^m(a_i^*a_i-a_ia_i^*)\|<1$. In particular
    we have
    \begin{equation}\label{one}
      \sum_{i=1}^m a_ia_i^*\le-1+t_0+\sum a_i^*a_i.
    \end{equation}
    Let $k=\|\sum_{i=1}^m a_i^*a_i\|$ and $a_{m+1}=(k-\sum_{i=1}^m
    a_i^*a_i)^{1/2}$. Then we have
    \begin{equation*}
      \sum_{i=1}^{m+1}a_i^*a_i=k
    \end{equation*}
    but on the other hand, by \eqref{one} we have also
    \begin{equation*}
      \sum_{i=1}^{m+1}a_ia_i^*\le -1+t_0+k.
    \end{equation*}
    The required properties are now fullfilled with $n=m+1$ and
    $b_i=(1/\sqrt{k})a_i$. 

    $(2)\Rightarrow(3)$. Suppose that $b_1,b_2,\ldots,b_n$ are as in
    $(3)$. Define $\phi(a)=\sum b_iab_i^*$. Then $\phi$ is a bounded
    positive map on $A$ with norm $\|\Phi\|=\|\sum b_ib_i^*\|<1$. It
    follows that $Id_A-\Phi$ is invertible in the Banach algebra
    $\mathcal{B}(A)$ of bounded maps on $A$. Let
    \begin{equation*}
      \Psi=(Id_A-\Phi)^{-1}.
    \end{equation*}
    Note that $\Psi=\sum_{i=0}^\infty\Phi^i$, thus $\Psi$ is positive too.
    By definition of $\Psi$, for any $a\in A$ we have
    \begin{equation*}
      a=(Id_A-\Phi)(\Psi(a))=\Psi(a)-\sum_{i=1}^n b_i\Psi(a)b_i^*=
      \sum_{i=1}^n[b_i^*,b_i\Psi(a)],
    \end{equation*}
    so $a$ is a finite sum of almost $n$ commutators
    If moreover $a$ is a positive element in $A$ then
    \begin{equation*}
      a=(Id_A-\Phi)(\Psi(a))=\Psi(a)-\sum_{i=1}^n b_i\Psi(a)b_i^*=
      \sum_{i=1}^n[\Psi(a)^{1/2}b_i^*,b_i\Psi(a)^{1/2}]
    \end{equation*}
    so $a$ is a finite sum of at most $n$ self-adjoint commutators.
  \end{proof}

\end{section}

\begin{section}{Questions}
  For a infinite $C^*$-algebra $A$ (in the sense of this paper) let $n(A)$
  be the least positive integer such that any element of $A$ is a sum of
  almost $n(A)$ commutators. In all the example we know, we have $n(A)=2$.
  We believe that it's unlikely to be always the case.
  
  In~\cite{haagerup} it was shown that, if $A$ is an unital exact
  $C^*$-algebra, then there exist an integer $m$ such that
  $\mathbb{M}_m(A)$ is properly infinite. It follows that
  $n(\mathbb{M}_m(A))=2$. Then a simple computation shows that
  $n(A)\le2m^2$. It would be interesting to answer the inverse problem,
  that is: assuming $n(A)$ is known, estimate the least positive integer
  such that $n(\mathbb{M}_n(A))=2$.
\end{section}

\bibliographystyle{amsplain} \bibliography{commutators} 
\end{document}